\documentclass{amsart}
\usepackage{amsmath} 
\usepackage{amssymb}
\usepackage{mathrsfs}

\newtheorem{theorem}{Theorem}[section] 
\newtheorem{claim}[theorem]{Claim}

\theoremstyle{definition}
\newtheorem{definition}[theorem]{Definition}

\theoremstyle{remark}

\newcommand{\gr}{{\rm gr}}

\newcommand{\rest}{{\restriction}}

\newcommand{\Wilog}{{\rm Without loss of generality}}

\newcommand{\then}{{\underline{then}}}

\newcommand{\mn}{{\medskip\noindent}}
\newcommand{\sn}{{\smallskip\noindent}}

\newcommand{\cA}{{\mathscr A}}

\newcommand{\gd}{{\mathfrak d\/}}

\newcommand{\bbN}{{\mathbb N}}

\newcommand{\frr}{{\mathfrak r}} 
\newcommand{\cU}{{\mathscr U}}

\newcommand{\cf}{{\rm cf}}

\newcount\skewfactor
\def\mathunderaccent#1#2 {\let\theaccent#1\skewfactor#2
\mathpalette\putaccentunder}
\def\putaccentunder#1#2{\oalign{$#1#2$\crcr\hidewidth
\vbox to.2ex{\hbox{$#1\skew\skewfactor\theaccent{}$}\vss}\hidewidth}}

\newenvironment{PROOF}[2][\proofname.]
   {\begin{proof}[#1]\renewcommand{\qedsymbol}{\textsquare$_{\rm #2}$}}
   {\end{proof}}

\begin{document}

\title {On reaping number having countable cofinality}
\author {Saharon Shelah}
\address{Einstein Institute of Mathematics\\
Edmond J. Safra Campus, Givat Ram\\
The Hebrew University of Jerusalem\\
Jerusalem, 91904, Israel\\
 and \\
 Department of Mathematics\\
 Hill Center - Busch Campus \\ 
 Rutgers, The State University of New Jersey \\
 110 Frelinghuysen Road \\
 Piscataway, NJ 08854-8019 USA}
\email{shelah@math.huji.ac.il}
\urladdr{http://shelah.logic.at}
\thanks{The author thanks Alice Leonhardt for the beautiful typing.
  First typed July 18, 2013. Paper E76}




\date{December 9, 2013 }

\begin{abstract}
We prove that if the bounding number $(\gd)$ is bigger than the
reaping number $(\frr)$, then the latter one has uncountable cofinality.
\end{abstract}

\maketitle
\numberwithin{equation}{section}
\setcounter{section}{-1}
\newpage

\section {Introduction}

\noindent
Recall
\begin{definition}
\label{z2}
Let $\gr$ the reaping number be the minimal cardinality of a set $\cA
\subseteq [\bbN]^{\aleph_0}$ such that for no set $X \subseteq N$ do
we have $A \in \cA \Rightarrow |A \cap X| = \aleph_0 = |A \backslash X|$.
\end{definition}
\newpage

\section {The proof}

\begin{claim}  
\label{a2}
If $\gd > \frr$ then $\cf(\frr) > \aleph_0$.
\end{claim}

\begin{PROOF}{\ref{a2}}  
Toward contradiction assume $\lambda = \frr = \sum\limits_{n}
\lambda_n$ where $\bigwedge\limits_{n} \lambda_n < \lambda_{n+1}$ (and
$\gd > \lambda$).  Let $\cA = \{A_\alpha:\alpha < \lambda\} \subseteq
[\omega]^{\aleph_0}$ witness $\frr \nleq \lambda$, i.e. $\cA$ is a
splitting family.  \Wilog \, $\alpha < \lambda \Rightarrow 
A_{2 \alpha +1} = \bbN \backslash A_{2\alpha}$; 
let $\cA_n = \{A_\alpha:\alpha < \lambda_n\}$.  Now 
we choose $\bar A^*_n = \langle
A^*_\eta:\eta \in {}^n 2\rangle$ by induction on $n$ such that
\mn
\begin{enumerate}
\item[$(*)_1$]  $(a) \quad \bar A^*_n$ is a partition of $\bbN$ to
infinite sets so $A_{<>} = \bbN$
\sn
\item[${{}}$]  $(b) \quad$ if $n=m+1$ and $\eta \in {}^n 2$ then
  $A^*_\eta \subseteq A^*_{\eta \rest m}$
\sn
\item[${{}}$]  $(c) \quad$ if $n=m+1,\eta \in {}^m 2$ and $\cU_\eta =
  \{\alpha < \lambda:A_\alpha \cap A^*_\eta$ infinite$\}$ \then \,

\hskip25pt $A^*_{\eta \char 94 <1>}$ (equivalently $A^*_{\eta \char 94 <0>}$)
  divide $A_\alpha \cap A^*_\eta$ to two 

\hskip25pt infinite sets when $\alpha \in \lambda_n \cap \cU_\eta$.
\end{enumerate}
\mn
This is possible as $\lambda_n < \frr$.  Having chosen $\bar A^*_n$ for
every $n$, for every $\eta \in {}^\omega 2$ and increasing $f:\omega
\rightarrow \omega$ we define $B_{\eta,f} = \cup\{A^*_{\eta \rest n}
\cap [0,f(n)):n <\omega\}$.

Now
\mn
\begin{enumerate}
\item[$(*)_2$]  if $\alpha < \lambda_n$ then there are $\eta_0,\eta_1
  \in {}^n 2$ such that

$\eta_1 \trianglelefteq \nu \in {}^{\omega >}2 \Rightarrow |A^*_\nu
  \cap A_{2 \alpha}| = \aleph_0$

$\eta_0 \trianglelefteq \nu \in {}^{\omega >}2 \Rightarrow |A^*_\nu
  \backslash A_{2 \alpha}| = \aleph_0$, equivalently $|A^*_\nu \cap
  A_{2 \alpha +1}| = \aleph_0$.
\end{enumerate}
\mn
[Why?  As $\bar A^*_n$ is a partition of $\bbN$ for some
  $\eta_0,\eta_1 \in {}^n 2$ we have $|A^*_{\eta_0} \cap A_{2
    \alpha}| = \aleph_0 = |A^*_{\eta_1} \cap A_{2 \alpha +1}|$ and
  then use the inductive choice of $\bar A^*_{n+1},\bar A^*_{n+2}
...$.  Note that really we have the partition of ${}^k 2$ to three 
$\langle v_{\alpha,k,\iota}:\iota < 3\rangle$ where $v_{\alpha,k,0}
= \{\eta \in {}^k 2:A^*_\eta \cap A_\alpha$ is finite$\}$,
$v_{\alpha,k,1} = \{\eta \in {}^k 2:A^*_\eta \backslash A_\alpha$
is finite$\}$; and note that $\eta \in {}^n 2 \wedge n \le k \wedge \eta
\trianglelefteq \nu \in {}^k 2 \wedge \iota < 3 \Rightarrow (\nu \in
v_{\alpha,k,\iota} \Rightarrow \eta \in v_{\alpha,n,\iota})$.]
\mn
\begin{enumerate}
\item[$(*)_3$]  if $\alpha < \lambda_n$ and $\eta \in {}^\omega 2,
|A_\alpha \cap A^*_{\eta \rest n}| = \aleph_0$ \then \, for some $g =
g_{\eta,\alpha} \in {}^\omega \omega$ we have: if $f \in {}^\omega
\omega$ and $f \nleq^* g$ then $|A_\alpha \cap B_{\eta,f}| = \aleph_0$.
\end{enumerate}
\mn
[Why?  Think, really we can make $f$ depend just on $A_\alpha$.]
\mn
\begin{enumerate}
\item[$(*)_4$]  if $\eta_1 \ne \eta_2 \in {}^\omega 2$ an $f_1,f_2 \in
  {}^\omega \omega$ then $B_{\eta_1,f_1} \cap B_{\eta_2,f_2}$ is
  finite.
\end{enumerate}
\mn
[Why?  Think.]
\mn
\begin{enumerate}
\item[$(*)_5$]  choose a sequence $\langle \eta_i:i < \omega\rangle$
  of pairwise distinct members of ${}^\omega 2$ which is dense,
  moreover $(\forall \nu \in {}^{\omega >}2)(\exists i)(\nu
  \triangleleft \eta_{2i} \wedge \nu \triangleleft \eta_{2i+1})$
\sn
\item[$(*)_6$]  for each $i$ we can choose $f_i:\omega \rightarrow
  \omega$ such that if $\alpha < \lambda_n$ and $A^*_{\eta_i \rest n}
  \cap A_\alpha$ is infinite then $B_{\eta_i,f_i} \cap A_\alpha$ is
  infinite.
\end{enumerate}
\mn
[Why?  By $\gd > \frr = \lambda > \lambda_n$ recalling the definition
  of $\gd$ and by $(*)_3$.]

Now let

\[
B_i = B_{\eta_i,f_i} \backslash \cup\{B_{\eta_j,f_j}:j <i\}
\]

\[
B^*_0 = \bigcup\limits_{i} B_{2i}
\]

\[
B^*_1 = \bigcup\limits_{i} B_{2i+1}
\]

\mn
$B^*_0,B^*_1$ are disjoint infinite $\subseteq \bbN$ and
satisfies $\iota < 2 \wedge \alpha
< \lambda \Rightarrow B^*_\iota \cap A_\alpha$ is infinite.

Why?  Let $\alpha < \lambda$ so for some $n,\alpha < \lambda_n$ hence
by $(*)_2$ for some $\rho_0 \ne \rho_1 \in {}^n 2$ we have 
\mn
\begin{enumerate}
\item[$\odot_1$]  $\rho_0 \triangleleft \nu \in {}^{\omega >}2
  \Rightarrow A^*_\nu \cap A_{2 \alpha}$ is infinite
\sn
\item[$\odot_2$]  $\rho_1 \triangleleft \nu \in {}^{\omega >}2
\Rightarrow A^*_\nu \cap A_{2 \alpha +1}$ is infinite.
\end{enumerate}
\mn
By $(*)_5$ there are $i_0,i_1$ such that $i_\ell = \ell \mod 2$ and
$\rho_\ell \triangleleft \eta_{i_\ell}$ for $\ell=0,1$.  Also the set
$B_{\eta_{i_0},f_{i_0}}$ satisfies $B_{\eta_{i_0},f_{i_0}} \cap A_{2
  \alpha}$ is infinite, but $B_{\eta_{i_0},f_{i_0}} \cap
\bigcup\limits_{j < i_0} B_{\eta_j,f_j}$ is finite hence $B_{i_0} \cap
A_{2 \alpha}$ is infinite but $i_0$ is even hence $B^*_0 \cap A_{2
  \alpha}$ is infinite.

Simiarly using $i_1$ which is odd $B^*_1 \cap A_{2 \alpha +1} = B^*_1
\backslash A_{2 \alpha}$ is infinite.
\end{PROOF}

\bibliographystyle{alphacolon}
\bibliography{lista,listb,listx,listf,liste,listz}

\end{document}